\documentclass{article}
\usepackage{fedembis}

\begin{document}

\title{Non-well-founded trees in categories \footnote{Paper submitted for publication.}}

\author{ Benno van den Berg \footnote{\it vdberg@math.uu.nl.} \and Federico De Marchi \footnote{\it marchi@math.uu.nl.} \thanks{Research supported by NWO Grant n.~613.000.222.}  \\ \\
Department of Mathematics, University of Utrecht \\
P.O.~Box 80010, 3508 TA Utrecht, The Netherlands}
\date{September 2004}

\maketitle

\begin{abstract}
  Non-well-founded trees are used in mathematics and computer science,
  for modelling non-well-founded sets, as well as non-terminating
  processes or infinite data-structures. Categorically, they arise as
  final coalgebras for polynomial endofunctors, which we call
  M-types. In order to reason about trees, we need the notion of path,
  which can be formalised in the internal logic of any locally
  cartesian closed pretopos with a natural number object. In such
  categories, we derive existence results about M-types, leading to
  stability of locally cartesian closed pretoposes with a natural
  number object and M-types under slicing, formation of coalgebras
  (for a cartesian comonad), and sheaves for an internal site.
\end{abstract}


\section{Introduction}

The relevance of non-well-founded trees to mathematics and computer
science was made clear by the work of Peter Aczel \cite{aczel88}. He
used them in order to extend the set-theoretic universe by
non-well-founded sets. Traditionally, the Axiom of Foundation allows
sets to be represented only by well-founded trees, but Aczel's
Anti-Foundation Axiom extends this possibility to all non-well-founded
trees.

In computer science, non-well-foundedness of trees enables one to
describe circular (and, more generally, non-terminating)
phenomena. For this reason, they have been used in the theory of
concurrency and specification, as well as in the study of semantics
for programming languages with coinductive types
\cite{cirstea00a,courcelle83,jacobsrutten96,ruttenturi98,barwisemoss89}.

Categorically, non-well-founded trees over a given signature form what
we call its M-type, that is the final coalgebra for the polynomial
functor determined by the signature itself, much like well-founded
trees form its initial algebra, usually called the W-type of the
signature.

In this paper, we look at properties of categories with
M-types. Specifically, we prove that they are closed under slicing,
formation of coalgebras (for a cartesian comonad), and sheaves (for an
internal site).

These constructions have proved useful in topos theory, leading to the
formulation of various independence results
\cite{freyd80,tierney72}. Our motivation for studying them in this
setting, is the potential for application to the theory of
non-well-founded sets and to polymorphism of coinductive types.

The core of the proof, in each case, is to show existence of certain
final coalgebras. Generally, this is not an easy task. It is obvious
that some functors can not have a final coalgebra (for example the
powerset functor in any elementary topos, by Cantor's
argument). However, playing with the structure of our categories and
restricting our attention to particular classes of functors, we can
obtain partial results (see for instance
\cite{aczelmendler89,barr99,worrell99,rutten00}).

In the present work, we shall focus on locally cartesian closed
pretoposes. The internal logic of these categories forces us to think
of trees in a constructive (and predicative) way. For this purpose, we
choose to use the language of paths \cite{courcelle83}; to this end we
shall require the existence of a natural number object. In formalising
our arguments, we encountered a category of objects, which we call
proto-coalgebras. We have no notice of them being ever introduced
before.

As a byproduct of our work with proto-coalgebras and paths, we could
reformulate in our setting a result by Santocanale
\cite{santocanale03}.

We are grateful to Jaap van Oosten and Ieke Moerdijk for their useful
comments on an earlier version of this paper.

\section{Setting the scene}\label{sec:prelim}

Throughout the paper, \ct{C} will denote a locally cartesian closed
pretopos with a natural number object. To any morphism \func{f}{B}{A}
in \ct{C}, we can then associate a \emph{polynomial functor}
\func{P_f}{\ct{C}}{\ct{C}}, defined as
$$P_f(X) = \Sigma_{a \in A} X^{B_a}$$
or, more formally, as
$$P_f(X)=\Sigma_A(A\times X \labto{p_1} A)^{(B \labto{f} A)},$$
where the exponential is taken in the slice category $\ct{C}/A$. In
the internal language of \ct{C}, an element of $P_f(X)$ is a pair
$(a,t)$, where $a\in A$ and \func{t}{f^{-1}(a)}{X} assigns to each
element in the fibre over $a$ an element in $X$. Notice in particular
that there is an obvious map \func{\rho}{P_f(X)}{A} mapping the pair
$(a,t)$ to $a$. This makes $P_f(X)$ into an object over $A$, namely
the exponential of $p_1$ and $f$. Notice also that the action of $P_f$
on maps lives over $A$, i.e.\ for a map \func{\phi}{X}{Y} in \ct{C}
the following commutes:
\diaglab{Pfpresroot}{
 P_f(X) \ar[rr]^{P_f\phi} \ar[rd]_{\rho}&& P_f(Y) \ar[ld]^{\rho} \\
 {} & A. &
}

The \emph{W-type} associated to $f$ is defined to be the initial
algebra for the functor $P_{f}$. In \ct{Set}, this is the set of all
well-founded trees built over the signature determined by $f$, having
one term constructor of arity $f^{-1}(a)$ for each element $a\in
A$. Initiality allows for definition of maps on trees by structural
induction. In this same context, the final $P_f$-coalgebra turns out
to be the set of all trees (including the non-well-founded ones) built
over the same signature. We shall call a final $P_f$-coalgebra the
\emph{M-type} associated to $f$, and we shall denote it by
$$\func{\tau_f}{M_f}{P_f(M_f)}.$$
We shall often denote the inverse of $\tau_f$ by $\sup_f$ (or just
$\sup$, when $f$ is understood). This will map the pair $(a,t)$ to the
tree $\sup_at$, consisting of a root labelled by $a$ and having
children indexed by the elements of the fibre of $f$ over $a$, in such
a way that the subtree appended to the branch labelled by $b\in
f^{-1}(a)$ is $t(b)$. Mediating arrows to final coalgebras correspond
to functions to trees defined by coinduction (for a clear description,
look at \cite{jacobsrutten96}). We say that \ct{C} \emph{has M-types},
if terminal coalgebras exist for every polynomial functor. A
\emph{coinductive pretopos} will be a locally cartesian closed
pretopos with a natural number object which has M-types.

It is the purpose of this paper to prove the closure of coinductive
pretoposes under slicing, formation of coalgebras for a cartesian
comonad and formation of sheaves over an internal site. But, before
proceeding, we state some general properties of M-types which will be
needed later.

\subsection{Covariant character of M-types}\label{sec:covariant}

First of all let's note that, by Lambek's lemma \cite{lambek70}, the
coalgebra map $\tau_f$ is an isomorphism. In particular, this
determines a \emph{root} map
\diag{M_f \ar[r]^-{\tau_f} & P_f(M_f) \ar[r]^-{\rho} & A,}
which, by an abuse of notation, will be denoted again by $\rho$. We are
confident that this will generate no confusion.

Given a pullback diagram in \ct{C}
\diag{B' \ar[d]_{f'} \ar[r]^-{\beta}& B \ar[d]^f \\
      A' \ar[r]_-{\alpha} & A,
}
we can think of $\alpha$ as a morphism of signatures, since the fibre
over each $a'\in A'$ is isomorphic to the fibre over $\alpha(a')\in
A$. It is therefore reasonable to expect, in such a situation, an
induced morphism between $M_{f'}$ and $M_f$, when these
exist.

In fact, as already pointed out in \cite{moerdijkpalmgren00}, such a
pullback square induces a natural transformation
\func{\widetilde{\alpha}}{P_{f'}}{P_f} such that
\begin{labequation}{tilderoot}
  \rho\widetilde{\alpha}=\alpha\rho.
\end{labequation}%
Post-composition with $\widetilde{\alpha}$ turns any $P_{f'}$-coalgebra
into one for $P_f$.  In particular, this happens for $M_{f'}$, thus
inducing a unique coalgebra homomorphism as in
\diaglab{covmtype}{
 M_{f'} \ar[r]^-{\alpha_!} \ar[d]_{\tau_{f'}} & M_f \ar[dd]^{\tau_f}\\
 P_{f'}(M_{f'}) \ar[d]_{\widetilde{\alpha}} & \\
 P_f(M_{f'}) \ar[r]_-{P_f(\alpha_!)} & P_f(M_f).
}
Notice that, by \refdiag{Pfpresroot} and \refeq{tilderoot}, the
morphism $\alpha_!$ preserves the root maps.

\subsection{Working with paths}
\label{sec:paths}

The reader who is familiar with the language of trees will not find it
surprising that, in order to carry out some calculations, we introduce
the notion of {\em path} \cite{courcelle83}. The reason for doing so
is that paths allow us to identify properties of trees in a
predicative way. Making an essential use of the internal logic of a
locally cartesian closed pretopos \ct{C} with a natural number object
\NN, we can define paths not just for the W-type and the M-type of a
polynomial functor, but for {\em any} coalgebra.

Assume we are given a $P_f$-coalgebra
$$ X \labto{\gamma} P_f X. $$
A finite sequence of odd length $< x_0, b_0, x_1, b_1, \ldots, x_n >$
is called a \emph{(finite) path} in $(X,\gamma)$, if every $x_i$ is in
$X$, every $b_i$ is in $B$ and for every $i\lt n$ we have
\begin{labequation}{relinpath}
  x_{i+1} = \gamma(x_i)(b_i).
\end{labequation}%
More precisely, if $\gamma(x_i)=(a_i,t_i)$, then we are asking that
$x_{i+1}=t_i(b_i)$. In the particular case when $X$ is the final
coalgebra $M_f$, a path $< m_0, b_0, \ldots, m_n >$ in this sense
coincides precisely with a path in the usual sense in the
non-well-founded tree $m_0$. We will therefore say that such a path
{\em lies} in $m_0$, and by extension, a path $<x_0,b_0,\dots,x_n>$
lies in $x_0\in X$ for any coalgebra $(X,\gamma)$. All paths in a
coalgebra $(X,\gamma)$ are collected in the subobject
\diag{ Paths(\gamma) \ar@{ >->}[r] & (X+B+1)^{\NN}.}

Any morphism of coalgebras \func{\alpha}{(X,\gamma)}{(Y,\delta)}
induces a morphism
\begin{labequation}{pathcoalghom}
  \func{\alpha^*}{Paths(\gamma)}{Paths(\delta)}
\end{labequation}%
between the respective objects of paths. A path $<x_0,b_0,\ldots,x_n>$
is sent by $\alpha^*$ to
$<\alpha(x_0),b_0,\ldots,\alpha(x_n)>$. Furthermore, given a path
$\tau =<y_0,b_0,\ldots,y_n>$ in $Y$ and an $x_0$ such that
$\alpha(x_0) = y_0$, there is a unique path $\sigma$ starting with
$x_0$ such that $\alpha^*(\sigma) = \tau$. (Proof: define $x_{i+1}$
inductively for every $i \lt n$ using \refeq{relinpath} and put
$\sigma = <x_0,b_0,\ldots,x_n >$.)

In fact, in order to introduce the concept of path, we need even less
than a coalgebra: it is sufficient to have a common environment in
which to read equation \refeq{relinpath}. Given a map \func{f}{B}{A}
in \ct{C}, consider the category $P_f-\ct{prtclg}$ of
$P_f$-\emph{proto-coalgebras}. Its objects are pairs of maps
\diaglab{protocol}{
  (\gamma,m) = X \ar[r]^-{\gamma} & Y & \ar@{ >->}[l]_-{m} P_f(X),
}
where $m$ is monic. An arrow between $(\gamma, m)$ and $(\gamma',m')$
is a pair of maps $(\alpha, \beta)$ making the following commute:
\diag{
  X \ar[r]^-{\gamma} \ar[d]_{\alpha} & Y \ar[d]_{\beta}
    & \ar@{ >->}[l]_-{m} P_f(X) \ar[d]^{P_f(\alpha)} \\
  X'\ar[r]_-{\gamma'} & Y' & \ar@{ >->}[l]^-{m'} P_f(X').
}
Notice that there is an obvious inclusion functor
\begin{labequation}{icoalgprotocol}
  \func{I}{\coalg{P_f}}{P_f-\ct{prtclg}},
\end{labequation}%
mapping a coalgebra \func{\gamma}{X}{P_f(X)} to the pair
$(\gamma,\id_{P_fX})$.

For a proto-coalgebra as in \refdiag{protocol}, one can introduce the
notion of a path in the following way. We shall call an element $x\in
X$ \emph{branching} if $\gamma(x)$ lies in the image of $m$. Then, we
call a sequence of odd length $\sigma= <x_0, b_0, x_1, b_1, \ldots,
x_n>$ a \emph{path} if it satisfies the properties:
\begin{enumerate}
  \item $x_i \in X$ is branching for all $i \lt n$; \\
  \item $b_i \in B_{a_i}$ for all $i \lt n$; \\
  \item $t_i(b_i) = x_{i+1}$ for all $i \lt n$;
\end{enumerate}
where $(a_i, t_i)$ is the (unique) element in $P_fX$ such that
$\gamma(x_i) = m(a_i, t_i)$. An element $x \in X$ is called
\emph{coherent}, if all paths starting with $x$ end with a branching
element. So, all coherent elements are automatically branching, and
their children, identified through $m$, are themselves coherent.

The internal language of \ct{C} makes it possible to identify the
object of coherent elements in any proto-coalgebra. This can be shown
to be the intersection of the chain $X_n$ of subobjects of $X$ defined
inductively as follows: $X_0=X$; $X_1$ is the pullback
\diag{
  X_1 \ar@{ >->}[r]^-{m_0} \ar[d] & X \ar[d]^{\gamma} \\
  P_f(X) \ar@{ >->}[r]_-m & Y;
}
given a subobject $\xymatrix@1{m_{n-1}:X_n\,\ar@{ >->}[r] & X_{n-1}}$,
$X_{n+1}$ is the pullback
\diag{
  X_{n+1}\ar@{ >->}[rrr]^-{m_n}\ar[d] &&&
    X_n\ar[d]^{\gamma m_0\dots m_{n-1}} \\
  P_f(X_n) \ar@{ >->}[rrr]_-{m P_f(m_0) \dots P_f(m_{n-1})} &&& Y.
}

Intuitively, a branching element is one whose image under $\gamma$ is
in $P_f(X)$, hence it has some \emph{branches} determined by the
fibres of $f$. However, nothing ensures that the children nodes will
be themselves branching. The subobject $X_n$ above is the object of
those elements in $X$ for which we can follow a path for at least
$n$-many steps, that is, we can find branching children up to $n$
consecutive generations. The intersection $Coh(\gamma)$ of the $X_n$
is then the object of coherent elements, for which we know that all
their children nodes are still coherent. Hence, we can
define on $Coh(\gamma)$ a $P_f$-coalgebra structure. In fact, this is
the biggest coalgebra which we can embed in $(\gamma,m)$, i.e.\ a
coreflection of the latter for the inclusion functor $I$ of
\refeq{icoalgprotocol}.

\begin{propo}{adjprotocol}
 The assignment $\xymatrix@1{(\gamma,m) \ar@{|->}[r] & Coh(\gamma)}$
 mapping any $P_f$-proto-coalgebra to the object of coherent elements
 in it, determines a right adjoint $Coh$ to the functor
 \func{I}{\coalg{P_f}}{P_f-\ct{prtclg}}.
\end{propo}

\begin{proof}
  Consider a proto-coalgebra
  \diag{ X \ar[r]^-{\gamma} & Y & \ar@{ >->}[l]_-{m} P_f(X), }
  and build the object $Coh(\gamma)$ of coherent elements in $X$.
  Because any coherent element $x\in Coh(\gamma)$ is also branching,
  we can find a (necessarily unique) pair $(a,t)$ such that
  $\gamma(x)=m(a,t)$. By defining $\chi(x)=(a,t)$ we equip
  $Coh(\gamma)$ with a $P_f$-coalgebra structure (notice that, $x$
  being coherent, so are the elements in the image of $t$). The
  coalgebra $(Coh(\gamma),\chi)$ clearly fits in a commutative diagram
  \diag{
    Coh(\gamma) \ar@{ >->}[rr]^-{i}\ar[d]_{\chi} && X\ar[d]^{\gamma}\\
    P_f(Coh(\gamma))\ar@{ >->}[r]_-{P_fi} & P_f(X)\ar@{ >->}[r]_-{m} & Y.
  } 

  Let now $(X',\chi')$ be any other $P_f$-coalgebra. Then, given a
  coalgebra morphism
  \diag{
    X' \ar[r]^-{\phi} \ar[d]_{\chi'} & Coh(\gamma) \ar[d]^{\chi} \\
    P_f(X') \ar[r]_-{P_f\phi} & P_f(Coh(\gamma)),
  }
  the pair $(i\phi, mP_f(i\phi))$ clearly determines a proto-coalgebra
  morphism from $I(X',\chi')$ to $(\gamma,m)$. Conversely, any
  proto-coalgebra morphism
  \diag{
    X' \ar[r]^-{\chi'} \ar[d]_{\alpha} & P_f(X') \ar[d]_{\beta}
      & \ar@{=}[l] P_f(X') \ar[d]^{P_f(\alpha)} \\
    X\ar[r]_-{\gamma} & Y & \ar@{ >->}[l]^-{m} P_f(X)
  }
  has the property that $\alpha(x')$ is branching for any $x'\in X'$.
  Using an opportune extension to proto-coalgebras of the morphism
  $\alpha^*$ described in \refeq{pathcoalghom} above, one can then
  easily check that elements in the image of $\alpha$ are
  coherent. Hence, $\alpha$ factors through the object $Coh(\gamma)$,
  inducing a coalgebra morphism from $(X',\chi')$ to
  $(Coh(\gamma),\chi)$.

  It is now easy to check that the two constructions are mutually
  inverse, thereby describing the desired adjunction.
\end{proof}

A particular subcategory of proto-coalgebras arises when we have
another endofunctor $F$ on \ct{C} and an injective natural
transformation $\xymatrix@1{m:P_f\,\ar@{ >->}[r]&F}$. In this case,
any $F$-coalgebra \func{\chi}{X}{FX} has an obvious
$P_f$-proto-coalgebra associated to it, namely $(\chi,m_X)$. The
assignment $(X,\chi)\mapsto (\chi,m_X)$ determines a functor
\func{\widehat{m}}{\coalg{F}}{P_f-\ct{prtclg}}, which is clearly
faithful.

\begin{propo}{restradjprotocol}
  The adjunction $I \ladj Coh$ of \refprop{adjprotocol} restricts to
  an adjunction $m_*\ladj Coh\;\widehat{m}$, where
  \func{m_*}{\coalg{P_f}}{\coalg{F}} takes \func{\chi}{X}{P_fX} to
  $(X,m_X\chi)$.
\end{propo}

\begin{proof}
  Consider a $P_f$-coalgebra $(Z,\gamma)$ and an $F$-coalgebra
  $(X,\chi)$. Then, a simple diagram chase, using the naturality of
  $m$, shows that $F$-coalgebra morphisms from $m_*(Z,\gamma)$ to
  $(X,\chi)$ correspond bijectively to morphisms of proto-coalgebras
  from $I(Z,\gamma)$ to $\widehat{m}(X,\chi)$, hence to
  $P_f$-coalgebra homomorphisms from $(Z,\gamma)$ to
  $Coh(\widehat{m}(X,\chi))$, by \refprop{adjprotocol}.
\end{proof}

\subsection{Existence of M-types}

The crucial point of the proof that coinductive pretoposes are closed
under the various constructions we are going to consider, will always
be that of showing existence of M-types. The machinery to do so will
be set up in this section.

Traditionally, one can recover non-well-founded trees from
well-founded ones, whenever the signature has one specified
constant. In fact, the constant allows for the definition of
truncation functions, which cut a tree at a certain depth and replace
all the term constructors at that level by that specified
constant. The way to recover non-well-founded trees is then to
consider sequences of trees $(t_n)_{n\gt 0}$ such that each $t_n$ is
the truncation at depth $n$ of $t_m$ for all $m\gt n$. Each such
sequence is viewed as the sequence of approximations of a given tree.

In our context, we call a map \func{f}{B}{A} \emph{pointed}, when the
signature it represents has a specified constant symbol, i.e.\ if
there exists a global element \func{\bot}{1}{A} such that the
following is a pullback:
\diag{
  0 \ar[r] \ar[d] & B \ar[d]^f \\
  1 \ar[r]_{\bot} & A.
}

The next two statements make clear that, instead of starting with
well-founded trees, i.e.\ with the W-type for $f$, we can build these
approximations from any fixed point of $P_f$.

\begin{lemm}{Mtypforpointedf}
  If for some pointed $f$ in \ct{C}, $P_f$ has a fixed point (that is,
  a (co)algebra for which the structure map is an isomorphism), then
  it also has a final coalgebra.
\end{lemm}
\begin{proof}
  Assume $X$ is an algebra whose structure map \func{\sup}{P_fX}{X} is
  an isomorphism. Observe, first of all, that $X$ has a global element
  \begin{equation}
    \func{\bot}{1}{X},
  \end{equation}
  namely $\sup_{\bot}(t)$, where $\bot$ is the point of $f$ and $t$ is
  the unique map from $B_{\bot}=0$ to $X$.

  Define, by induction, the following truncation functions
  \func{tr_n}{X}{X}:
  \begin{eqnarray*}
    tr_0 & = & \bot \\
    tr_{n+1} & = & \sup\circ P_f(tr_n)\circ\textstyle{\sup^{-1}}
  \end{eqnarray*}

  Using these maps, we can define the object $M$, consisting of
  all sequences $(\alpha_n \in X)_{n \gt 0}$ with the property:
  $$\alpha_n = tr_n(\alpha_m) \textrm{ for all } n \lt m.$$

  On $M$, we define a $P_f$-coalgebra structure \func{\tau}{M}{P_f M}
  as follows. Given a sequence $\alpha = (\alpha_n)\in M$, observe
  that $\rho(\alpha_n)$ is independent of $n$ and is some element $a
  \in A$. Hence, each $\alpha_n$ is of the form $\sup_a(t_n)$ for some
  \func{t_n}{B_a}{X}, and we define \func{t}{B_a}{M} by putting
  $t(b)_n = t_{n+1}(b)$ for every $b \in B_a$; then $\tau(\alpha) =
  (a, t)$. Thus, $M$ has the structure of a $P_f$-coalgebra, and we
  claim it is the terminal one.

  To show this, given another coalgebra \func{\chi}{Y}{P_fY}, we wish
  to define a map of coalgebras \func{\widehat{p}}{Y}{M}. This means
  defining maps \func{\widehat{p}_n}{Y}{X} for every $n\gt 0$, with
  the property that $\widehat{p}_n = tr_n \widehat{p}_m$ for all $n\lt
  m$. Intuitively, $\widehat{p}_n$ maps a state of $Y$ to its
  ``unfolding up to level $n$'', which we can mimic in $X$. Formally,
  they are defined inductively by
  \begin{eqnarray*}
    \widehat{p}_0 & = & \bot \\
    \widehat{p}_{n+1} & = & \sup\circ P_f(\widehat{p}_n)\circ\chi.
  \end{eqnarray*}

  It is now easy to show, by induction on $n$, that $\widehat{p}_n=
  tr_n\widehat{p}_m$ for all $m\gt n$. For $n=0$, both sides of the
  equation become the constant map $\bot$. Supposing the equation
  holds for a fixed $n$ and any $m\gt n$, then for $n+1$ and any $m\gt
  n$ we have 
  \begin{eqnarray*}
    \widehat{p}_{n+1} & = & \sup\,P_f(\widehat{p}_{n})\,\chi  \\
      & = & \sup\,P_f(tr_n\,\widehat{p}_m)\,\chi \\
      & = & \sup P_f(tr_n)\,\textstyle{\sup^{-1}}\,\sup\,P_f(\widehat{p}_m)\,            \chi \\
      & = & tr_{n+1}\,\widehat{p}_{m+1}.
  \end{eqnarray*}
  We leave to the reader the verification that $\widehat{p}$ is the
  unique $P_f$-coalgebra morphism from $X$ to $M$.
\end{proof}

\begin{theo}{fixptimpmtype}
  If fixed points exist in \ct{C} for all $P_f$ (with $f$ pointed),
  then \ct{C} has M-types .
\end{theo}
\begin{proof}
  Let \func{f}{B}{A} be a map. We freely add a point to the
  signature represented by $f$, by considering the composite
  \diaglab{pointedf}{
    f_{\bot}:B \ar[r]^-f & A \ar@{ >->}[r]^-{i} & A + 1
  }
  (with the point \func{j = \bot}{1}{A + 1}). Notice that the obvious
  pullback
  \diag{
    B\ar[r]^-{\id} \ar[d]_f & B \ar[d]^{f_{\bot}} \\
    A \ar@{ >->}[r]_-i & A+1
  }
  determines, by \refdiag{covmtype}, a (monic) natural transformation
  \func{i_!}{P_f}{P_{f_{\bot}}}; hence, by \refprop{restradjprotocol},
  the functor \func{(i_!)_*}{\coalg{P_f}}{\coalg{P_{f_{\bot}}}} has a
  right adjoint. Now observe that $P_{f_{\bot}}$ has a fixed point, by
  assumption, hence a final coalgebra by \reflemm{Mtypforpointedf}.
  This will be preserved by the right adjoint of $(i_!)_*$, hence
  $P_f$ has a final coalgebra.
\end{proof}

This proof gives a categorical counterpart of the standard
set-theoretic construction: add a dummy constant to the signature,
build infinite trees by sequences of approximations, then select the
actual M-type by taking those infinite trees which involve only term
constructors from the original signature. This last passage is
performed by the coreflection functor of \refprop{restradjprotocol},
since branching elements are trees in the M-type of $f_{\bot}$ whose
root is not $\bot$, and coherent ones are trees with no occurrence of
$\bot$ at any point.

From this last theorem, we readily deduce the following result, first
pointed out to us by Abbott, Altenkirch and Ghani (this was, in fact,
one of the motivations for looking at this topic). Here, we call
\emph{predicative topos} a locally cartesian closed pretopos with
W-types.

\begin{coro}{WtypesMtypes} \cite{abbottaltenkirchghani04}
  Every predicative topos is a coinductive pretopos.
\end{coro}
\begin{proof}
  Every predicative topos \ct{C} has a natural number object, namely
  the W-type associated to the left inclusion
  \func{\inl}{1}{1+1}. Since the W-type associated to a (pointed) map
  $f$ is a fixed point for $P_f$, \ct{C} also has all M-types by the
  previous theorem.
\end{proof}

Because, as pointed out in \cite{moerdijkpalmgren00}, an elementary
topos with a natural number object has all W-types, so is a
predicative topos, this also establishes (in a very roundabout
fashion):

\begin{corollary} \cite{fiveauthors01}
  Any elementary topos with a natural number object has M-types.
\end{corollary}

Together, these corollaries provide us with a substantial class of
examples of coinductive pretoposes. It is still an open question
whether there is any relevant examples of a coinductive pretopos that
is not a predicative topos.

Although \reftheo{fixptimpmtype} is clearly helpful in proving that
certain categories have M-types, it is even more so, when combined
with the following observation. Given any monic morphism
\func{\alpha}{P_fX}{X}, we can think of it as identifying branching
elements in $X$, giving at the same time an explicit description of
their children. We can therefore form the chain $(X_n)_{n\gt 0}$ of
``elements in $X$ with $n$ generations of children''. The intersection
of the $X_n$ (the object of coherent elements in $X$) will then be a
$P_f$-coalgebra; in fact, a fixed point for $P_f$.

\begin{propo}{injalgimpfixpt}
  Let \ct{C} be a locally cartesian closed pretopos with a natural
  number object, and \func{f}{B}{A} a map in it. Then, any prefixed
  point \func{\gamma}{P_fX}{X} (i.e.\ an algebra whose structure map
  is monic) has a subalgebra that is a fixed point.
\end{propo}
\begin{proof}
  Any prefixed point \func{\alpha}{P_fX}{X} can be seen as a
  $P_f$-proto-coalgebra
  \diag{ 
    X\ar[r]^-{\id} & X & \ar@{ >->}[l]_-{\alpha} P_fX.
  }
  Its coreflection $Coh(\id,\alpha)$, as defined in
  \refprop{adjprotocol}, is a $P_f$-coalgebra \func{\gamma}{Y}{P_fY}
  (in fact, the largest) fitting in the following commutative square:
  \diag{
    Y \ar@{ >->}[r]^-i \ar[d]_{\gamma} & X \\
    P_fY \ar@{ >->}[r]_-{P_fi} & P_fX. \ar@{ >->}[u]_{\alpha}
  }
  Now, consider the image under \func{I}{\coalg{P_f}}{P_f-\ct{prtclg}}
  of the coalgebra \func{P_f(\gamma)}{P_fY}{P_f^2Y}. The morphism
  of proto-coalgebras
  \diag{
    P_fY \ar[r]^-{P_f\gamma} \ar[d]_{\alpha P_fi} & 
      P_f^2Y \ar[d]|{\alpha P_f(\alpha)P_f^2i} &
      \ar@{ >->}[l]_-{\id} P_f^2Y \ar[d]^{P_f(\alpha)P_f^2i} \\
    X \ar[r]_-{\id} & X & \ar@{ >->}[l]^-{\alpha} P_fX
  }
  transposes through $I\ladj Coh$ to a morphism
  \func{\phi}{(P_fY,P_f\gamma)}{(Y,\gamma)}, which is a right inverse
  of \func{\gamma}{(Y,\gamma)}{(P_fY,P_f\gamma)} by the universal
  property of $(Y,\gamma)$. Hence, we have
  $\gamma\phi=P_f(\phi\gamma)=\id$, proving that $\gamma$ and $\phi$
  are mutually inverse.
\end{proof}

Putting together \reftheo{fixptimpmtype} and \refprop{injalgimpfixpt},
we get at once the following:

\begin{coro}{injalgimpmtype}
  If in a locally cartesian closed pretopos \ct{C} with a natural
  number object we have prefixed points for every polynomial functor,
  then \ct{C} has M-types.
\end{coro}

As an application of the techniques in this section, we present the
following result, which is to be compared with the one by Santocanale
in \cite{santocanale03}. An immediate corollary of his Theorem 4.5
is that M-types exist in every locally cartesian closed pretopos
with with a natural number object, for maps of the form \func{f}{B}{A}
where $A$ is a finite sum of copies of $1$. Notice that such an object
$A$ has {\em decidable equality}, i.e.\ the diagonal
\func{\Delta}{A}{A\times A} has a complement in the subobject lattice
of $A \times A$. We extend the statement above to {\em all} maps whose
codomain has decidable equality.

\begin{propo}{decidableequality}
  Let \ct{C} be a locally cartesian closed pretopos with a natural
  number object, and \func{f}{B}{A} a morphism in \ct{C}. If $A$ has
  decidable equality, then the M-type for $f$ exists.
\end{propo}
\begin{proof}
  Without loss of generality, we may assume that $f$ is pointed; in
  fact, if we replace $A$ by $A_{\bot} = A + 1$ and $f$ by $f_{\bot}$
  as in \refdiag{pointedf}, then $A_{\bot}$ also has decidable
  equality, and the existence of an M-type for the composite
  $f_{\bot}$ implies that of an M-type for $f$ (see the proof of
  \reftheo{fixptimpmtype}). Then, by \refprop{injalgimpfixpt} and
  \reflemm{Mtypforpointedf}, it is enough to show that $P_f$ has a
  prefixed point.

  Let $S$ be the object of all finite sequences of the form
  $$ < a_0, b_0, a_1, b_1, \ldots , a_n > $$
  where $f(b_i) = a_i$ for all $i \lt n$. (Like paths in a coalgebra,
  this object $S$ can be constructed using the internal logic of
  \ct{C}.) Now, let $V$ be the object of all decidable subobjects of
  $S$ (these can be considered as functions $S \rTo 1 + 1$). Define
  the map \func{m}{P_fV}{V} taking a pair $(a, \func{t}{B_a}{V})$ to
  the subobject $P$ of $S$ defined by the following clauses:
  \begin{enumerate}
    \item $< a_0 > \in P$ iff $a_0 = a$.
    \item $<a_0,b_0>*\sigma\in P$ iff $a_0=a$ and $\sigma\in t(b_0)$.
  \end{enumerate}
  (Here, $*$ is the symbol for concatenation.) $P$ is obviously
  decidable, so $m$ is well-defined. To see that it is monic, suppose
  $P=m(a,t)$ and $P'=m(a',t')$ are equal. Then,
  $$ <a> \in P \Longrightarrow <a> \in P' \Longrightarrow a = a', $$
  and, for every $b \in B_a$ and $\sigma \in S$,
  \begin{eqnarray*}
    \sigma \in t(b) & \Longleftrightarrow & <a, b> * \sigma \in P \\
                    & \Longleftrightarrow & <a, b> * \sigma \in P' \\
                    & \Longleftrightarrow & \sigma \in t'(b),
  \end{eqnarray*}
  so $t = t'$ and $m$ is monic. Hence, $(V, m)$ is a prefixed point for
  $P_f$ and we are finished.
\end{proof}

\begin{rema}{decidableelements}
  To obtain the M-type for $f$ from $V$, one should first deduce a
  fixed point $V'$ from it, as in \refcoro{injalgimpmtype}. This means
  selecting the coherent elements of $V$, and these turn out to be those
  decidable subobjects $P$ of $S$ satisfying the following properties:
  \begin{enumerate}
    \item $<a> \in P$ for a unique $a \in A$;
    \item if $<a_0,b_0,\ldots,a_n>\in P$, then there exists a unique
      $a_{n+1}$ for any $b_n\in B_{a_n}$ such that
      $<a_0,b_0,\ldots,a_n,b_n,a_{n+1}>\in P$.
  \end{enumerate}
  Now, we should turn this fixed point into the M-type for $f$ (as in
  \reflemm{Mtypforpointedf}), but this step is redundant, since our
  choice of $V$ is such that $V'$ already is the desired M-type.
\end{rema}

\section{M-types and slicing}

We start by considering preservation of the coinductive pretopos
structure under slicing. Let $I$ be an object in a locally cartesian
closed pretopos with a natural number object \ct{C}. Then, it is
well-known that the slice category $\ct{C}/I$ has again the same
structure, and the reindexing functor
\func{x^*}{\ct{C}/I}{\ct{C}/J} for any map \func{x}{J}{I} in \ct{C}
preserves it. So, we can focus on showing the existence of M-types in
$\ct{C}/I$ and their preservation under reindexing.

We shall first concentrate on the existence, proving a ``local
existence'' result, from which we derive a global statement. Then, we
shall investigate the action of the reindexing functors.

Let us consider a map
\diaglab{mapslice}{
  B \ar[rr]^-{f} \ar[rd]_{\beta} && A \ar[ld]^{\alpha} \\
  & I
}
in $\ct{C}/I$. We shall denote by $P_f$ the polynomial functor
determined by $f$ (or, more precisely, by $\Sigma f$) in \ct{C}, and
by $P^I_f$ the polynomial endofunctor determined in $\ct{C}/I$.
The functor \func{P_f}{\ct{C}}{\ct{C}} can be extended to a functor
\func{P_f}{\ct{C}}{\ct{C}/I}; in fact, $P_fX$ lives over $A$ via the
root map, and the composite \func{\alpha\rho}{P_fX}{I} defines the
desired extension.

\begin{lemm}{PIfsubfunctor}
  There is an injective natural transformation of endofunctors on
  $\ct{C}/I$
  $$\func{i}{P^I_f}{P_f\Sigma_I.}$$
\end{lemm}

\begin{proof}
  Fix an object \func{\xi}{X}{I} in $\ct{C}/I$. Then, unfolding the
  definitions, we have
  \begin{eqnarray*}
    \Sigma_IP^I_f(X\labto{\xi}I) & = & \Sigma_I\Sigma_{\alpha}
      (A\times_I X \labto{\alpha^*\xi} A)^{(B\labto{f}A)} \\
    & = & \Sigma_A(A\times_I X \labto{\alpha^*\xi}A)^{(B\labto{f}A)}
  \end{eqnarray*}
  and
  $$P_f(\Sigma_I(X\labto{\xi}I)) = P_f(X) =
    \Sigma_A(A\times X\labto{p_1} A)^{(B\labto{f}A)},$$
  where \func{\alpha}{A}{I} makes $A$ into an object over $I$. The
  pullback $A\times_IX$ fits into the equaliser diagram
  \diag{
    A\times_IX\ar[r]^-{\psi}&
    A\times X \ar@<.5ex>[r]^-{\alpha p_1} \ar@<-.5ex>[r]_-{\xi p_2}
    & I,
  }
  hence $\psi$ is a map over $A$, since $p_1\psi=\alpha^*\xi$.
  Moreover, $\psi$ is monic in $\ct{C}/A$, and so is $\Sigma_A\psi^f$,
  thus giving the desired monomorphism, which is easily seen to live
  over $I$. Naturality is also readily checked.
\end{proof}
\begin{rema}{PIfequaliser}
  In the internal language of \ct{C}, the map $i_{(X,\xi)}$ realises
  $\Sigma_IP^I_f(X,\xi)$ as the subobject of $P_f(X)$ consisting of
  those pairs $(a,\func{t}{f^{-1}(a)}{X})$ such that for all $b\in
  f^{-1}(a)$ it holds $\alpha(a)=\xi t(b)$. In other words, $i$ is the
  equaliser
  {\xymatrixrowsep={1em}\diag{
    \Sigma_IP^I_f(X,\xi)\ar[r]^-i
      & P_fX \ar[rr]^-{P_f\xi} \ar[rd]_{\rho} & & P_fI, \\
      & & A \ar[ru]_k &
  }}
  where $k$ maps $a\in A$ to $(a,\lambda b.\alpha(a))$.

  If we think of the elements of $A$ as term constructors, this is
  just saying that $P^I_f(X,\xi)$ applies each $a$ just to elements of
  $X$ which sit in the same fibre over $I$.
\end{rema}

Using the $i$ of \reflemm{PIfsubfunctor}, we can build an M-type for
$f$ in $\ct{C}/I$, whenever $M_f$ exists in \ct{C}.

\begin{theo}{localmtypeslice}
  Let \ct{C} be a locally cartesian closed pretopos with a natural
  number object and $I$ an object in \ct{C}. Consider a map
  \func{f}{B}{A} over $I$, such that the functor
  \func{P_f}{\ct{C}}{\ct{C}} has a final coalgebra. Then, $f$ has an
  M-type in $\ct{C}/I$.
\end{theo}
\begin{proof}
  Let \func{\tau_f}{M_f}{P_fM_f} be the M-type associated to $f$ in
  \ct{C}. $M_f$ can be considered as an object over $I$, by taking the
  composite $\mu$ of the root map \func{\rho}{M_f}{A} with the map
  \func{\alpha}{A}{I}, and $(M_f,\tau_f)$ then becomes the final
  $P_f\Sigma_I$-coalgebra, as one can easily check. The adjunction
  determined by the natural transformation
  \func{i}{P^I_f}{P_f\Sigma_I} as in \refprop{restradjprotocol} takes
  the final $P_f\Sigma_I$-coalgebra $(M_f,\tau_f)$ to its coreflection
  $M^I_f$, and because right adjoints preserve limits, this is the
  final $P^I_f$-coalgebra.
\end{proof}

\begin{remark}
  \rm The injective natural transformation $i$ of
  \reflemm{PIfsubfunctor} identifies branching elements in
  $P_f\Sigma_I$ as those obtained by applying a term constructor in
  $A$ to elements living in its same fibre over $I$, as discussed in
  \refrema{PIfequaliser}.

  The coreflection process used to build $M^I_f$ out of the M-type
  $(M_f,\tau_f)$, helps understanding which elements of the latter do
  actually belong to the former. Trees in $M^I_f$ are coherent for the
  notion of branching determined by $P^I_f$, hence, not only the
  children of the root node live in its same fibre over $I$, but all
  the children of the children do too, and so on for any node in the
  tree. In other words, $M^I_f$ consists of those trees in $M_f$ all
  nodes of which live in the same fibre over $I$. As such, the object
  $M^I_f$ can also be described as the equaliser
  {\xymatrixrowsep={1em}\diag{
    M^I_f \ar@{ >->}[r] & 
      M_f \ar[rr]^-{<\id,\alpha>_!} \ar[rd]_{<\id,\alpha\rho>} & &
      M_{f\times\!I}, \\
      & & M_f\!\times I\! \ar[ru]_{\chi} &
  }}
  where $\chi$ is the map coinductively defined as
  $$\textstyle{\chi(\sup_at,i)=\sup_{(a,i)}(\chi<t,i>).}$$
\end{remark}

As an immediate consequence of \reftheo{localmtypeslice}, we get the
following:

\begin{corollary}
  For any given object $I$ of a coinductive pretopos \ct{C}, the slice
  category $\ct{C}/I$ is again a coinductive pretopos.
\end{corollary}

\begin{remark} \rm
  This last result could have also been proved directly by combining
  \refcoro{injalgimpmtype} and \reflemm{PIfsubfunctor}. However, the
  proof of \reftheo{localmtypeslice} shows that the construction is
  actually simpler. More specifically, notice that, in this case, we
  obtain the M-type for a map $f$ directly after the coreflection, and
  we do not need to add any dummy variable, nor to build sequences of
  approximations.
\end{remark}

We now look at the reindexing functors. Suppose we are given a
morphism \func{x}{J}{I} in \ct{C}, and consider the induced functor
\func{x^*}{\ct{C}/I}{\ct{C}/J}, with left adjoint $\Sigma_x$ (we shall
denote by $\eta$ the unit of this adjunction). Recall that the
reindexing functors, in this case, do always preserve the cartesian
closed structure of each slice.

Now, consider a morphism $f$ as in \refdiag{mapslice}, inducing a map
\func{x^*\!f}{x^*B}{x^*A} in $\ct{C}/J$ and a pullback square
\diag{
  x^*A \ar[r]^-{\alpha^*x} \ar[d]_{x^*\alpha} & A \ar[d]^{\alpha} \\
  J \ar[r]_-x & I.}

By the Beck-Chevalley property for this square, and the fact that
$x^*$ preserves exponentials, we then have, for any object
$(X,\xi)$ in $\ct{C}/I$,
\begin{eqnarray*}
  x^*P_f^I(X,\xi)
    & = & x^*\Sigma_{\alpha}(X\times_IA\rTo A)^{(B\labto{f}A)} \\
    & = & \Sigma_{x^*\alpha}
      (\alpha^*x)^*(X\times_IA\rTo A)^{(B\labto{f}A)} \\
    & \cong & \Sigma_{x^*\alpha}
      (\alpha^*x)^*(X\times_IA\rTo A)^{(\alpha^*x)^*(B\labto{f}A)} \\
    & = & \Sigma_{x^*\alpha}
      (x^*X\times_Jx^*A\rTo x^*A)^{(x^*B\labto{x^*\!f}x^*A)} \\
    & = & P_{x^*\!f}^J(x^*X,x^*\xi).
\end{eqnarray*}
Let us write \func{\theta}{P_{x^*\!f}^Jx^*X}{x^*P_f^IX} for the
isomorphism, and $\widehat{\theta}$ for its transpose under the
adjunction $\Sigma_x\ladj x^*$.

Given any coalgebra \func{\gamma}{X}{P_f^I X} in $\ct{C}/I$, we can
map it to a coalgebra
$$x^*X\labto{x^*\gamma}x^*P_f^IX\labto{\theta^{-1}}P_{x^*\!f}^Jx^*X.$$
This determines a functor \func{G}{\coalg{P_f^I}}{\coalg{P_{x^*\!f}^J}}.

Now, given a coalgebra \func{\delta}{Y}{P_{x^*\!f}^JY} in $\ct{C}/J$, a
morphism \func{\phi}{(Y,\delta)}{G(X,\gamma)} determines a commutative
square
\diag{
  Y \ar[rr]^-{\phi} \ar[d]_{\delta} && x^{*}X \ar[d]^{x^{*}\gamma} \\
  P_{x^{*}\!f}^JY \ar[r]_-{P_{x^{*}\!f}^J\phi} &
    P_{x^{*}\!f}^Jx^{*}X \ar[r]_-{\theta} & x^{*}P_f^IX,
}
which transposes under $\Sigma_x\ladj x^*$ to the square
$$\xymatrix@C=3.5em{
  \Sigma_x \ar[d]_{\Sigma_x\delta} \ar[rr]^-{\widehat{\phi}} && 
    X \ar[d]^{\gamma} \\
  \Sigma_xP_{x^{*}f}^JY \ar[r]_-{\Sigma_xP_{x^{*}f}^J\phi} &
    \Sigma_xP_{x^{*}f}^Jx^{*}X \ar[r]_-{\widehat{\theta}} & P_f^IX.
}$$
By simple diagram chasing, one shows that
\begin{eqnarray*}
  \gamma\widehat{\phi}
    & = & \widehat{\theta}\,\Sigma_xP_{x^*\!f}^J(\phi)\,\Sigma_x(\delta) \\
    & = & \widehat{\theta}\,\Sigma_xP_{x^*\!f}^Jx^*(\widehat{\phi})\,
          \Sigma_xP_{x^*\!f}^J(\eta)\,\Sigma_x(\delta) \\
    & = & P_f^I(\widehat{\phi})\,\widehat{\theta}\,
          \Sigma_xP_{x^*\!f}^J(\eta)\,\Sigma_x(\delta);
\end{eqnarray*}
hence, the functor assigning to $(Y,\delta)$ the composite
$$\xymatrix@1{
  \Sigma_xY \ar[r]^-{\Sigma_x\delta} & 
    \Sigma_xP_{x^*\!f}^JY \ar[rr]^-{\Sigma_xP_{x^*\!f}^J\eta} &&
    \Sigma_xP_{x^*\!f}^Jx^*\Sigma_xY \ar[r]^-{\widehat{\theta}} &
    \Sigma_xY
}$$
is a left adjoint to $G$.

In particular, this implies that the image of the M-type
$(M_f^I,\tau_f)$ for $f$ in $\ct{C}/I$, along the right adjoint $G$ is
the final $P_{x^*\!f}^J$-coalgebra in $\ct{C}/J$, i.e.\
$G(M_f^I)=x^*M_f^I$ is the M-type for $x^*\!f$. This proves the following.

\begin{theo}{reindexing}
  For any map \func{x}{J}{I} in a coinductive pretopos \ct{C}, the
  reindexing functor \func{x^*}{\ct{C}/I}{\ct{C}/J} preserves the
  structure of a coinductive pretopos.
\end{theo}

\section{M-types and coalgebras}

In this section, we turn our attention to the construction of
categories of coalgebras for a cartesian comonad $(G,\epsilon,\delta)$
(by \emph{cartesian} we mean that $G$ preserves finite limits). As for
the slicing case, we already know that most of the structure of a
coinductive pretopos is preserved by taking coalgebras for $G$:

\begin{theorem}
  If \ct{C} is a locally cartesian closed pretopos with natural number
  object, then so is $\ct{C}_G$ for a cartesian comonad $G = (G,
  \epsilon, \delta)$ on \ct{C}.
\end{theorem}

\begin{proof}
  Theorem 4.2.1 on page 173 of \cite{johnstone02} gives us that
  $\ct{C}_{G}$ is cartesian, in fact locally cartesian closed, and
  that it has a natural number object. The two additional requirements
  of having finite disjoint sums and being exact are easily verified,
  since the forgetful functor \func{U}{\ct{C}_G}{\ct{C}} creates
  finite limits.
\end{proof}

Given a morphism $f$ of $G$-coalgebras, this determines a
polynomial functor \func{P_f}{\ct{C}_G}{\ct{C}_G}, while its
underlying map $Uf$ determines the endofunctor $P_{Uf}$ on \ct{C}. The
two are related as follows:

\begin{propo}{Pfsubcoalg}
  Let \func{f}{(B,\beta)}{(A,\alpha)} be a map of $G$-coalgebras.
  Then, there is an injective natural transformation
  \diag{
    \ct{C} \ar@{}[rd]|{\stackrel{i}{\rightarrowtail}}
       \ar[d]_-G \ar[r]^-{P_{U\!f}} & \ct{C} \ar[d]^G \\
    \ct{C}_G \ar[r]_-{P_f} & \ct{C}_G,
  }
  whose mate under the adjunction $U\ladj G$ we shall denote by
  \diaglab{matetransf}{
    j:UP_f\ar[r] & P_{Uf}U: \ct{C}_G \ar[r] & \ct{C}.
  }
\end{propo}
In order to prove this result, we first need to introduce a couple of
lemmas (the proof of the second can also be found in \cite{johnstone02}).
\begin{lemm}{expcoalg}
  Let \ct{C} be a cartesian closed category and $G$ be a cartesian
  comonad on it. Then, for any object $X$ in \ct{C} and any coalgebra
  $Y$, the exponential $(GX)^Y$ exists and is equal to $G(X^{UY})$.
\end{lemm}
\begin{proof}
  The forgetful functor \func{U}{\ct{C}_G}{\ct{C}} is left adjoint to
  the cofree coalgebra functor \func{G}{\ct{C}}{\ct{C}_G}. Moreover,
  $U$ preserves products because it creates them; hence, the following
  bijective correspondence:
  \begin{displaymath}\renewcommand{\arraystretch}{1.2}\begin{array}{rcl}
    A & \rTo & G(X^{UY}) \\ \hline
    UA & \rTo & X^{UY} \\ \hline
    UA \times UY & \rTo & X \\ \hline
    U(A \times Y) & \rTo & X \\ \hline
    A \times Y & \rTo & GX.
  \end{array}\end{displaymath}
\end{proof}
\begin{lemm}{equivslices}
  Every $G$-coalgebra $(A,\alpha)$ in $\ct{C}_G$, determines a
  cartesian comonad $G'$ on $\ct{C}/A$ and an isomorphism between
  $\ct{C}_{G}/(A, \alpha)$ and $(\ct{C}/A)_{G'}$ which respects the
  forgetful functors \func{U}{\ct{C}_{G}/(A, \alpha)}{\ct{C}/A} and
  \func{U'}{(\ct{C}/A)_{G'}}{\ct{C}/A}.
\end{lemm}
\begin{proof}
  The comonad $G'$ is computed on an object $t: X \rTo A$ in
  $\ct{C}/A$ by taking the following pullback:
  \diaglab{defofgprime}{
    G'X \ar[d]_{G't} \ar@{ >->}[r] & GX \ar[d]^{Gt} \\
    A \ar@{ >->}[r]_{\alpha}       & GA.
  }
  Then, there is a natural one-to-one correspondence between the
  dotted arrows over $A$ and $GA$ in
  \diag{
    X \ar @/_/ [ddr]_t \ar @/^/@{.>}[drr]^s \ar@{.>}[dr]|r \\
    & G'X \ar[d]^{G't} \ar@{ >->}[r] & GX \ar[d]^{Gt} \\
    & A \ar@{ >->}[r]_-{\alpha}       & GA.
  }
  It is now easy to show that if $r$ is a $G'$-coalgebra, then the
  corresponding $s$ is a $G$-coalgebra, and vice versa. Hence, the
  isomorphism. Moreover,
  $$U(\func{t}{(X,s)}{(A,\alpha)})=t=U'(t,r);$$
  therefore, the action of the forgetful functors is respected.

  Note that both horizontal arrows in the pullback
  \refdiag{defofgprime} are monic, because $\epsilon_A$ is a
  retraction of the $G$-coalgebra $\alpha$.
\end{proof}

{\setlength{\parindent}{0pt}
{\bf Proof of \refprop{Pfsubcoalg}.}
  Recall that $P_f(GX)$ is the source of the exponential $(A \times GX
  \rTo A)^f$ in the category $\ct{C}_G/(A,\alpha)$. However, this can
  be more easily computed in the category of $G'$-coalgebras. First of
  all, since $G$ preserves products, the object $A\times GX \rTo A$
  corresponds, through the isomorphism of \reflemm{equivslices}, to
  $G'(A\times X \rTo A)$. Then, using \reflemm{expcoalg}, the
  exponential takes the form $G'((A\times X\rTo A)^{Uf})$, and fits in
  the following pullback square, which is an instance of
  \refdiag{defofgprime}:
  \diag{
    G'((A\times X\to A)^{U\!f}) \ar[d] \ar@{ >->}[r]^{i_X}
      & G((A\times X\to A)^{U\!f}) \ar[d] \\
    A \ar@{ >->}[r]_{\alpha} & GA.
  }
  Now notice that the top-right entry of the diagram is exactly
  $GP_{U\!f}(X)$, hence the map $i$ therein defines the $X$-th
  component of the desired natural transformation.
\mbox{} \hfill $\Box$ \mbox{}\\}

We are now ready to formulate a local existence result for M-types in
categories of coalgebras.

\begin{theo}{localmtypecoalg}
 Let \func{f}{(B,\beta)}{(A,\alpha)} be a map of $G$-coalgebras. If
 the underlying map $Uf$ has an M-type in \ct{C}, then the functor
 \func{P_f}{\ct{C}_G}{\ct{C}_G} has a final coalgebra in $\ct{C}_G$.
\end{theo}
\begin{proof}
  The natural transformation $i$ of \refprop{Pfsubcoalg} allows us to
  turn any $P_{Uf}$-coalgebra into a $P_f$-proto-coalgebra. In
  particular, for the M-type \func{\tau}{M=M_{Uf}}{P_{Uf}M} in \ct{C},
  we get the proto-coalgebra
  \diag{
    GM \ar[r]^-{G\tau} & GP_{Uf}M & \ar@{ >->}[l]_-{i_M} P_fGM,
  }
  whose coreflection $Coh(M)=Coh(G\tau,i_M)$ is final in \coalg{P_f}.
  To see this, consider another coalgebra $(X,\gamma)$ (therefore, $X$
  is a $G$-coalgebra, and \func{\gamma}{X}{P_fX} is a homomorphism of
  $G$-coalgebras). To give a morphism of $P_f$-coalgebras from
  $(X,\gamma)$ to $Coh(M)$ is the same, through $I\ladj Coh$, as
  giving a map \func{\psi}{X}{GM} in $\ct{C}_G$ which is a morphism of
  $P_f$-proto-coalgebras, i.e. that makes the following commute:
  \diag{
    X \ar[rr]^-{\gamma} \ar[d]_{\psi} && P_fX \ar[d]^{P_f\psi} \\
    GM \ar[r]_-{G\tau} & GP_{Uf}M & \ar@{ >->}[l]^-{i_M} P_fGM.
  }
  This transposes, through $U\ladj G$, to the following diagram in
  \ct{C}, where $j$ is the natural transformation defined in
  \refdiag{matetransf}:
  \diag{
    UX \ar[r]^-{U\gamma} \ar[d]_{\widehat{\psi}} &
      UP_fX \ar[r]^-{j_X} & P_{Uf}UX \ar[d]^{P_{Uf}\widehat{\psi}}\\
    M \ar[rr]_-{\tau} && P_{Uf}M.
  }
  But finality of $M$ implies that there is precisely one such
  $\widehat{\psi}$ for any coalgebra $(X,\gamma)$, hence finality is
  proved.
\end{proof}
\begin{coro}{Mtypescoalg}
  If \ct{C} is a coinductive pretopos and $G=(G,\epsilon,\delta)$ is a
  cartesian comonad on \ct{C}, then the category $\ct{C}_G$ of
  (Eilenberg-Moore) coalgebras for $G$ is again a coinductive pretopos.
\end{coro}

\begin{remark}
  \rm Notice that \refcoro{Mtypescoalg} could also be deduced by
  \refcoro{injalgimpmtype}, in conjunction with
  \refprop{Pfsubcoalg}. However, analogously to what happens in the
  slicing case, \reftheo{localmtypecoalg} shows that we do not need to
  perform the whole construction, since the coreflection step gives
  directly the final coalgebra.
\end{remark}

\begin{remark}
  \rm In particular, this result shows stability of coinductive
  pretoposes under the glueing construction, since this is a special
  case of taking coalgebras for a cartesian comonad (see
  \cite{wraith74}).
\end{remark}

\section{M-types and sheaves}

In this section, we turn our attention to formation of sheaves for an
internal site in a coinductive pretopos. Our aim is to show that the
resulting category is again a coinductive pretopos. The structure of
the proof mimics the one for W-types in \cite{moerdijkpalmgren00}, but
we now make use of the finality of M-types and, implicitly, of the
language of paths.

We start by considering an internal category \CC\ in a coinductive
pretopos \ct{E}, with object of objects $\CC_0$ and domain and
codomain morphisms indicated by \func{d_0,d_1}{\CC_1}{\CC_0},
respectively. We shall assume \CC\ to have pullbacks and to come
equipped with a collection $T$ of covering families, satisfying the
axioms of a site \cite[Def.\,2 p.111]{maclanemoerdijk94}.

Notice that the category of presheaves \pshct{\CC} is then the
category of coalgebras for a cartesian comonad on the slice category
$\ct{E}/\CC_0$ (see for instance \cite[Example
A.4.2.4~($b$)]{johnstone02}). By the results of the previous sections,
we get at once

\begin{proposition}
  The presheaf category \pshct{\CC} is a coinductive pretopos.
\end{proposition}

Before proceeding to the sheaf case, we find it useful to have a
closer look at the construction of M-types for categories of
presheaves. This will provide a concrete description of the M-types,
and at the same time set some notation about polynomial functors in
presheaf (and sheaf) categories.

\subsection{M-types in presheaves}

First of all, we need to introduce the functor
\func{|\cdot|}{\pshct{\CC}}{\ct{E}} which takes a presheaf
\psh{A} to its ``underlying object'' $|\psh{A}|=\{ (a,C) \mid
a\in\psh{A}(C)\}$. This is just the composite of the forgetful functor
\func{U}{\pshct{\CC}}{\ct{E}/\CC_0} with
\func{\Sigma_{\CC_0}}{\ct{E}/\CC_0}{\ct{E}}.

Let \func{f}{\psh{B}}{\psh{A}} be a morphism of presheaves. Then, the
``fibre'' $\psh{B}_a$ of $f$ over $a\in\psh{A}(C)$ for an object $C$
in \CC\ is a presheaf, whose action on $D$ is described in the
internal language of \ct{E} as
$$\psh{B}_a(D)=\{(\beta,b)\mid\func{\beta}{D}{C},a\cdot\beta=f(b)\}$$
and restriction along a morphism \func{\delta}{D'}{D} is defined
as
$$(\beta,b)\cdot\delta=(\beta\delta,b\cdot\delta).$$

The image $P_f(X)$ of a presheaf $X$ is defined on an object $C$ of
\CC\ as
\begin{labequation}{Pfpresh}
  P_f(X)(C)=\{(a,t)\mid a\in\psh{A}(C),\,\func{t}{\psh{B}_a}{X}\},
\end{labequation}%
where $t$ is a presheaf morphism. Restriction along a morphism
\func{\alpha}{C'}{C} is defined as
$$(a,t)\cdot\alpha = (a\cdot\alpha,\,\alpha^*(t)),$$
where $\alpha^*(t)(\beta,b)=t(\alpha\beta,b)$.

But the presheaf morphism $f$ also induces a map
$$\func{f'}{\Sigma_{(a,C)\in|\psh{A}|}|\psh{B}_a|}{|\psh{A}|}$$
whose fibre over $(a,C)$ is precisely $|\psh{B}_a|$. This induces a
polynomial endofunctor on \ct{E}, which acts on an object $Y$ as
$$P_{f'}Y=
  \{((a,C),t)\mid a\in\psh{A}(C),\,\func{t}{|\psh{B}_a|}{Y}\}.$$
Notice that we can always think of $P_{f'}Y$ as the underlying object
of a presheaf, whose section over $C$ are those elements of $P_{f'}Y$
of the form $((a,C),t)$. Given a morphism \func{\alpha}{C'}{C} for
which $a\cdot\alpha=a'$, this induces a map
\func{\widetilde{\alpha}}{|\psh{B}_{a'}|}{|\psh{B}_a|}. Restriction
along $\alpha$ is then defined by $((a,C),t)\cdot\alpha=
((a',C'),t\widetilde{\alpha})$.

In particular, $f'$ induces an endofunctor
\func{F}{\pshct{\CC}}{\pshct{\CC}}, which verifies the equation
\begin{labequation}{Pfprpsh}
  |FX|=P_{f'}|X|=
  \{((a,C),t)\mid a\in\psh{A}(C),\,\func{t}{|\psh{B}_a|}{|X|}\}
\end{labequation}%
for any presheaf $X$ on \CC. Because each morphism of presheaves
$\psh{B}_a\rTo X$ induces a map on the underlying objects, it is clear
by \refeq{Pfpresh} and \refeq{Pfprpsh} that $P_f$ is a subfunctor of
$F$. We shall denote by $m$ the corresponding injective natural
transformation.

By \refprop{restradjprotocol}, we can then determine an M-type for $f$
by giving a final coalgebra for the functor $F$ and taking its
coreflection in \coalg{P_f}. The latter is readily obtained from the
M-type $(M_{f'},\tau_{f'})$ for $f'$ in \ct{E}. In fact, as we saw
above, the object $P_{f'}M_{f'}$ has a presheaf structure, and we can
translate it on $M_{f'}$ via the isomorphism $\tau_{f'}$, which then
becomes a presheaf morphism. To see that it is the final
$F$-coalgebra, it is now enough to observe that, given any other
$F$-coalgebra \func{\gamma}{X}{FX}, finality of $M_{f'}$ determines a
unique morphism of $P_{f'}$-coalgebras
\diag{
  |X| \ar@{..>}[r] \ar[d]_{|\gamma|} & M_{f'} \ar[d]^{\tau_{f'}} \\
  P_{f'}|X|=|FX| \ar@{..>}[r] & P_{f'}M_{f'},
}
which is easily seen to be a morphism of presheaves.

We have just given a proof of the following local existence result:

\begin{theorem}
  Consider a map \func{f}{\psh{B}}{\psh{A}} in \pshct{\CC}. If the
  induced map $f'$ has an M-type in \ct{E}, then $f$ has an M-type in
  \pshct{\CC}.
\end{theorem}

\begin{remark}
  \rm The natural transformation $m$, in this case, selects as
  branching only those trees $\sup_{(a,C)}t$ in $M_{f'}$ for which the
  map \func{t}{|\psh{B}_a|}{M_{f'}} is the underlying map of a
  presheaf morphism. The reader who is familiar with the language of
  \cite{moerdijkpalmgren00}, will recognise that trees coherent for
  $m$ are those called \emph{natural} there.
\end{remark}

\subsection{M-types in sheaves}

Once again, it is already known that the category of sheaves for an
internal site in a locally cartesian closed pretopos with a natural
number object has again the same structure, so, as in the previous
cases, all we need to show is that \shct{\CC} has M-types. Analogously
to what is done in \cite{moerdijkpalmgren00}, we get the result by
showing the following.
\begin{proposition}
  Let \func{f}{\psh{B}}{\psh{A}} be a morphism of sheaves, and
  $\psh{M}$ the corresponding M-type in \pshct{\CC}. Then, \psh{M} is
  a sheaf.
\end{proposition}
\begin{proof}
  The proof will work as follows. First, we shall define a presheaf
  $\psh{M}^+$, which is the separated presheaf associated to \psh{M}
  (see \cite{maclanemoerdijk94}). Then, we will show that $\psh{M}^+$
  is a $P_f$-coalgebra with mediating morphism
  \func{\widehat{m}}{\psh{M}^+}{\psh{M}}; moreover, the unit
  \func{\eta}{\psh{M}}{\psh{M}^+} is a coalgebra homomorphism. This
  will imply that $\widehat{m}\eta=\id_{\psh{M}}$, hence $\eta$ is
  mono, i.e.\ \psh{M} is separated. Existence of the glueing of a
  compatible family will be shown to be determined by $\widehat{m}$,
  and this will prove \psh{M} to be a sheaf.

  The presheaf $\psh{M}^+$ is defined on an object $C$ in \CC\ as
  $$\psh{M}^+(C) = \left\{(\{C_i\labto{c_i}C\},\,T_i\in\psh{M}(C_i))\:\Big|
    \begin{array}{l}
      \{c_i\} \mbox{ is a covering family and} \\[-1.5ex]
      \mbox{the $T_i$'s are compatible}
    \end{array}
    \right\}/\sim$$
  where
  $$(\{C_i\labto{c_i}C\},T_i)\,\sim\,(\{C'_j\labto{c'_j}C\},S_j)$$
  if there is a common refinement of $\{c_i\}$ and $\{c'_j\}$ on which
  the $T_i$'s and the $S_j$'s agree. More precisely, we require the
  existence of a covering family $\{\func{c''_k}{C''_k}{C}\}$ such
  that every $c''_k$ factors as $c_i\sigma_{ki}$ and $c'_j\tau_{kj}$
  for opportune maps $\sigma_{ki}$ and $\tau_{kj}$, with
  $T_i\cdot\sigma_{ki}=S_j\cdot\tau_{kj}$. We shall indicate the
  equivalence class of a pair $(\{c_i\},T_i)$ by enclosing it in
  square brackets.

  By the stability axiom for a site, a given morphism
  \func{\beta}{D}{C} in \CC\ induces, via pullback, a restriction of a
  covering family $\{c_i\}$ over $C$ to a covering family $\{d_i\}$ as
  below:
  \diaglab{restrcompfam}{
    D_i \ar[r]^-{\beta_i} \ar[d]_{d_i} & C_i \ar[d]^{c_i} \\
    D \ar[r]_-{\beta} & C.
  }
  We can then define the restriction morphism of $\psh{M}^+$ by letting
  \begin{labequation}{restrM+}
    [\{c_i\},T_i]\cdot\beta = [\{d_i\},T_i\cdot\beta_i].
  \end{labequation}%
  We now need to show that $-\cdot\beta$ is a well-defined map and
  that it satisfies the functoriality properties, thus making
  $\psh{M}^+$ into a presheaf. We shall go through this calculation in
  some detail, in order to illustrate the kind of arguments
  used. Analogous forms of reasoning will occur several times in this
  proof. All the successive details will then be omitted.

  Consider the two related pairs
  $$(\{c_i\},T_i)\sim(\{c'_j\},S_j).$$
  It is clearly enough to consider the case where the covering family
  $\{c'_j\}$ is a refinement of the family $\{c_i\}$, say via maps
  $\sigma_{ji}$, in such a way that for each $i$ there is a $j$ with
  $T_i\cdot\sigma_{ji}=S_j$. Pulling back the families $\{c_i\}$ and
  $\{c'_j\}$ along $\beta$, we get covering families $\{d_i\}$ and
  $\{d'_j\}$ as below:
  \diag{
    D_i \ar[r]^-{\beta_i} \ar[d]_{d_i} & C_i \ar[d]^{c_i}
      && D'_j\ar[r]^-{\beta'_j} \ar[d]_{d'_j} & C'_j \ar[d]^{c'_j} \\
    D \ar[r]_-{\beta} & C && D \ar[r]_-{\beta} & C.
  }
  The maps $\sigma_{ji}$ pull back to maps $\tau_{ji}$ such that
  $d'_j=d_i\tau_{ji}$, hence $\{d'_j\}$ is a refinement of $\{d_i\}$.
  Moreover, we have
  $$(T_i\cdot\beta_i)\cdot\tau_{ji}=T_i\cdot(d_i\tau_{ji})
    =T_i\cdot(\sigma_{ji}\beta'_j)=(T_i\cdot\sigma_{ji})\cdot\beta'_j
    =S_j\cdot\beta'_j,$$
  whence $[\{d_i\},T_i\cdot\beta_i]=[\{d'_j\},S_j\cdot\beta'_j]$. This
  proves that restriction along $\beta$ is well-defined.

  Now we show that this restriction operation defines a presheaf
  structure on $\psh{M}^+$. Preservation of identities is easily
  verified. As for composition, consider an element $[\{c_i\},T_i]$ in
  $\psh{M}^+(C)$ and composable maps
  $$E\labto{\delta}D\labto{\beta}C.$$
  Then, we want to show that
  \begin{labequation}{dotbetapresh}
    (\{e_i\},T_i\cdot\gamma_i)\sim(\{e'_i\},T_i\cdot\beta_i\delta_i),
  \end{labequation}%
  where the maps $e_i$, $\gamma_i$, $e'_i$, $\beta_i$ and $\delta_i$
  arise as shown by the pullback diagrams below:
  \diag{
    E_i \ar[r]^-{\gamma_i} \ar[d]_{e_i} & C_i \ar[d]^{c_i} &&
      E'_i \ar[r]^-{\delta_i} \ar[d]_{e'_i} & D_i \ar[r]^-{\beta_i}
      \ar[d]^{d_i} & C_i \ar[d]^{c_i}  \\
    E \ar[r]_{\beta\delta} & C &&
      E \ar[r]_{\delta} & D \ar[r]_{\beta} & C.
  }
  By pullback-pasting, there are obvious induced isomorphisms
  \func{\sigma_i}{E'_i}{E_i} such that
  $\gamma_i\sigma_i=\beta_i\delta_i$ and
  $e_i\sigma_i=e'_i$. Therefore, the family $\{e'_i\}$ is a refinement
  of $\{e_i\}$; moreover,
  $$T_i\cdot\beta_i\delta_i=T_i\cdot\gamma_i\sigma_i,$$
  which proves \refeq{dotbetapresh}.

  So, $\psh{M}^+$ is a presheaf, and there is an obvious presheaf
  morphism \func{\eta}{\psh{M}}{\psh{M}^+}, whose component on an
  object $C$ maps a tree $T\in\psh{M}(C)$ to $[\{\id_C\},T]$. Equation
  \refeq{restrM+} immediately implies that $\eta$ is natural. We shall
  soon show that it is also a coalgebra morphism, but first we have to
  equip $\psh{M}^+$ with a $P_f$-coalgebra structure.

  For this, we need a presheaf morphism
  \func{m}{\psh{M}^+}{P_f\psh{M}^+}. This will associate to an element
  $F=[\{\func{c_i}{C_i}{C}\},T_i]$ in $\psh{M}^+(C)$ a pair $(a,t)$
  with $a$ in $\psh{A}(C)$ and \func{t}{\psh{B}_a}{\psh{M}^+} a
  presheaf morphism. Suppose for each $i$ that
  $T_i=\sup_{(a_i,C_i)}t_i$, with \func{t_i}{\psh{B}_{a_i}}{\psh{M}}.
  Then, since the $T_i$'s form a compatible family, so do the $a_i$'s,
  thus determining a unique amalgamation $a$ such that $a\cdot
  c_i=a_i$, because \psh{A} is a sheaf. In order to define the action
  of $t$, consider an element $(\beta,b)$ in $\psh{B}_a(D)$ (so that
  \func{\beta}{D}{C} and $a\cdot\beta=f(b)$). The compatible family
  $\{c_i\}$ restricts along $\beta$ to a compatible family $\{d_i\}$
  as in \refdiag{restrcompfam}, and we define
  $$ t(\beta,b)=[\{d_i\},t_i(\beta_i,b\cdot d_i)]. $$
  By a standard calculation (of the kind shown above), it is easy to
  show that $t$ defines a presheaf morphism and that the definition of
  both $a$ and $t$ does not depend on the choice of the representative
  for $F$. Hence, the morphism $m$ is well-defined, and again by
  similar arguments, we can then show that it defines a presheaf
  morphism.

  So, we now have a coalgebra $m$ on the presheaf $\psh{M}^+$, for
  which it is easily checked that the presheaf morphism $\eta$ defined
  above is in fact a coalgebra homomorphism. Moreover, there is a
  unique map of coalgebras \func{\widehat{m}}{\psh{M}^+}{\psh{M}},
  determined by the finality of \psh{M}. The composite
  $\widehat{m}\eta$ is then a coalgebra morphism from \psh{M} to
  itself, therefore it must be the identity. In particular, this
  implies that $\eta$ is monic, i.e.\ \psh{M} is separated.  To
  complete the proof of the statement, we now need to show that every
  compatible family has an amalgamation. In fact, it turns out that
  the glueing of a family $F$ is determined by its image under
  $\widehat{m}$.

  More precisely, given a covering family
  $\{\func{c_i}{C_i}{C}\}$ over $C$ and a matching family of trees
  $T_i\in\psh{M}(C_i)$, we want to show that
  $T=\widehat{m}([\{c_i\},\,T_i])$ is the amalgamation of the $T_i$'s;
  that is, $T\cdot c_i=T_i$ for all $i$. To this end, fix an index
  $i_0$, and restrict the covering family $\{c_i\}$ along the map
  $c_{i_0}$ as below:
  \diaglab{gluerestr}{
    D_i \ar[r]^-{c'_i} \ar[d]_{d_i} & C_i \ar[d]^{c_i} \\
    C_{i_0} \ar[r]_-{c_{i_0}} & C.
  }
  We then have that
  \begin{eqnarray*}
    \widehat{m}([\{c_i\},T_i])\cdot c_{i_0}
      & = & \widehat{m}([\{c_i\},T_i]\cdot c_{i_0}) \\
    & = & \widehat{m}([\{d_i\},T_i\cdot c'_i])
  \end{eqnarray*}
  But clearly, the covering family $\{d_i\}$ also forms a refinement
  of $\{\id_{C_{i_0}}\}$, and since the $T_i$'s are a matching family
  (hence, in particular, for the pullbacks in \refdiag{gluerestr} we
  have that $T_{i_0}\cdot d_i=T_i\cdot c'_i$), it follows at once that
  $$(\{d_i\},T_i\cdot c'_i)\sim(\{\id_{C_{i_0}}\},T_{i_0}),$$
  whence, because $\widehat{m}\eta=\id$,
  $$T\cdot c_{i_0} = \widehat{m}([\{c_i\},T_i])\cdot c_{i_0}
    = \widehat{m}([\{\id_{C_{i_0}}\},T_{i_0}])
    = \widehat{m}\eta(T_{i_0}) = T_{i_0}.$$
\end{proof}

\begin{remark} \label{exvreg}
  \rm Notice that this is the only proof we have given where we make
  use of the exactness of coinductive pretoposes. In any other
  argument regularity would have been enough.
\end{remark}

Piecing this result together with the remarks at the beginning of this
section, we get at once the following:

\begin{theorem}
  Let \ct{E} be a coinductive pretopos and $(\CC,T)$ an internal site in
  \ct{E}. Then, the category \shct{\CC,T} of sheaves over \CC\ in \ct{E}
  is again a coinductive pretopos.
\end{theorem}

\section{Concluding remarks}

The original motivation for writing this paper was to set the
categorical framework in which to work with non-well-founded
trees. The natural continuation of this research is two-folded, with
applications to set theory and type theory.

The reason to expect such relationships to exist is the
well-established correspondence between categories with W-types,
Martin-L\"of type theory, and Aczel's constructive set theory (the
formal system $CZF$) (see for example
\cite{moerdijkpalmgren00,moerdijkpalmgren02,aczel78}). We suggest that
this correspondence should have a counterpart in the non-well-founded
world, between coinductive pretoposes, Martin-L\"of type theory with
coinductive types, and Aczel's non-well-founded set theory
($CZF^-+AFA$). This would provide the ground for successive research
in various directions.

\refcoro{WtypesMtypes} reveals that there are links between the
well-founded and the non-well-founded world. In particular, it shows
that within Martin-L\"of type theory we can use well-founded types to
model certain coinductive types. An analogous connection arises from
the work of Lindstr\"om \cite{lindstrom89}, where it is shown how to
construct non-well-founded sets in Martin-L\"of type theory. Making
the correspondence above clear, should bring the two approaches
together.

Ultimately, we expect our result on sheaves to lead to
independence results in non-well-founded set theory, and to
Martin-L\"of type theory with coinductive types.

\end{document}